\documentclass[12pt]{article}

\usepackage[latin1]{inputenc}
\usepackage[french]{babel}

\usepackage{amsfonts}
\usepackage{amssymb}
\newtheorem{lemme}{Lemme}[section]
\newtheorem{proposition}[lemme]{Proposition}
\newtheorem{theoreme}[lemme]{Th\'eor\`eme}
\newtheorem{corollaire}[lemme]{Corollaire}
 
\newcommand\dem{\noindent{\it D\'emonstration.}\ }
\newcommand\findem{\hfill $\square$ }
\newcommand\hfld[2]{\smash{\mathop{\hbox to 12mm{\rightarrowfill}}
     \limits^{\scriptstyle#1}_{\scriptstyle#2}}}
\newcommand\hflg[2]{\smash{\mathop{\hbox to 6mm{\leftarrowfill}}
     \limits^{\scriptstyle#1}_{\scriptstyle#2}}}
\newcommand\vfld[2]{\llap{$\scriptstyle#1$}
     \left\downarrow\vbox to 4mm{}\right.\rlap{$\scriptstyle #2$}}
\newcommand\vflm[2]{\llap{$\scriptstyle#1$}\left\uparrow\vbox to
4mm{}\right.
     \rlap{$\scriptstyle #2$}}
\newcommand\hfl[1]{\,\smash{\mathop{\hbox to 6mm{\rightarrowfill}}
     \limits^{\scriptstyle#1}}\,}

\newcommand\rta{\rightarrow}

\newcommand\calW{{\mathcal W}}
\newcommand\calV{{\mathcal V}}
\newcommand\calL{{\mathcal L}}
\newcommand\calM{{\mathcal M}}

\newcommand\calO{{\mathcal O}}
\newcommand\RR{{\mathbb R}}
\newcommand\NN{{\mathbb N}}
\newcommand\ZZ{{\mathbb Z}}
\newcommand\QQ{{\mathbb Q}}
\newcommand\GG{{\mathbb G}}
\newcommand\bbA{{\mathbb A}}

\newcommand\rmR{{\rm R}}
\newcommand\rmH{{\rm H}}

\newcommand\Aut{{\rm Aut}}
\newcommand\Spec{{\rm Spec}}

\newcommand\smv{{\scriptscriptstyle\vee}}

\newcommand\Id{{\rm Id}}
\newcommand\KR{{\rm KR}}

\title{Alcoves et $p$-rang des vari\'et\'es ab\'eliennes}
\author{ B.C. Ng\^o, A. Genestier}
\date{\small d\'edi\'e au Professeur Ng\^o Huy C\^an}
\begin{document}
\maketitle

\begin{abstract}
On \'etudie la relation entre le $p$-rang des vari\'et\'es ab\'eliennes 
en caract\'eristique $p$ et la stratification de Kottwitz-Rapoport 
de la fibre sp\'eciale en $p$ de l'espace de module des vari\'et\'es 
ab\'eliennes principalement polaris\'ees avec structure de niveau 
de type Iwahori en $p$. En particulier, on d\'emontre la densit\'e 
du lieu ordinaire dans cette fibre sp\'eciale.

\bigskip
\centerline{\bf Alcoves and $p$-rank of abelian varieties.} 

\bigskip
We study the relation between the $p$-rank of abelian varieties 
in characteristic $p$ and the Kottwitz-Rapoport's stratification 
of the special fiber modulo $p$ of the moduli space of principally polarized 
abelian varieties with Iwahori type level structure on $p$.
In particular, the density of the ordinary locus in that special 
fiber is proved.
\end{abstract}

\section*{Introduction}

Dans \cite{deJong}, de Jong a \'etudi\'e la mauvaise r\'eduction  modulo
$p$ de l'espace de modules ${{\cal A}}$ des vari\'et\'es ab\'eliennes
de  dimension $n$, principalement polaris\'ees et avec structure de
niveau
$\Gamma_0(p)$ en $p$. En utilisant le th\'eor\`eme de
Grothendieck-Messing sur les d\'eformations de sch\'emas ab\'eliens, il
a ramen\'e l'\'etude des singularit\'es de cette mauvaise r\'eduction
\`a l'\'etude de celles d'un autre probl\`eme de modules $\calM$ d\'efini
p\^urement en termes d'alg\`ebre lin\'eaire. Il s'agit l\`a d'un cas
particulier des {\em mod\`eles locaux}  que Rapoport et Zink \cite{RZ}
ont associ\'e aux vari\'et\'es de Shimura de type PEL.

L'une des propri\'et\'es fondamentales de $\calM$ est que sa fibre 
sp\'eciale peut \^etre naturellement plong\'ee dans l'ind-sch\'ema de
drapeaux affines (au sens de \cite{BL}) du groupe symplectique ({\it cf.}
\cite{KR}, \cite{Goertz},
\cite{HN} et aussi \cite{Gen}). De ce fait, elle est  r\'eunion (finie,
indic\'ee par une partie finie $\KR$ du groupe de Weyl affine) d'orbites
sous l'action du sous-groupe d'Iwahori. Ceci induit sur la fibre
sp\'eciale
${{\cal A}}$ une certaine stratification. Le propos de cette note est de
d\'emontrer que cette stratification
est plus fine que celle d\'efinie par le
$p$-rang des vari\'et\'es ab\'eliennes. On obtiendra en particulier que
les strates ordinaires sont celles associ\'ees aux \'el\'ements de $\KR$
qui sont des translations. En conjonction avec un th\'eor\`eme de
Kottwitz-Rapoport, ceci d\'emontre donc que les strates ordinaires sont
denses dans la fibre sp\'eciale de ${{\cal A}}$, ce qui r\'epond \`a
une question pos\'ee par de Jong dans
\cite{deJong}. 

En ce qui concerne l'organisation du papier, les trois premi\`eres
sections o\`u nous rappelons le contexte du probl\`eme, sont connues des
experts. Seule la derni\`ere section est  donc originale --signalons
toutefois que
U. Goertz a aussi trouv\'e la formule pour le $p$-rang, alors que
circulait d\'ej\`a une premi\`ere version de cet article dans laquelle
il n'\'etait  question que du probl\`eme de densit\'e.

Nous exprimons notre gratitude \`a M.  Rapoport duquel  nous
avons beaucoup appris au sujet de la mauvaise r\'eduction des
vari\'et\'es de Shimura et qui par ses commentaires \cite{Rapoport} a
contribu\'e \`a am\'eliorer cet article. Nous remercions T.  Haines pour
avoir lu attentivement le manuscript et pour de nombreuses discussions
sur ce sujet. Nous remercions aussi R.  Kottwitz et J.  Tilouine pour
l'inter\^et qu'ils ont port\'e \`a ce travail.


\section{Rappels sur le mod\`ele local}

On fixe un entier $n\geq 1$, un nombre premier $p$ et un entier $N$
premier
\`a $p$. On note $T=\Spec(\ZZ_p)$, $\eta$ le point g\'en\'erique de $T$
et
$s$ son point ferm\'e. 

On consid\`ere le foncteur ${{\cal A}}$ qui associe \`a toute
$\ZZ_p$-alg\`ebre
$R$, l'ensemble des classes d'isomorphisme de 
$$A=\left(A_0\hfl\alpha A_1\hfl\alpha\cdots\hfl\alpha A_n,
\lambda_0,\lambda_n,\iota_N\right)$$ o\`u
\begin{itemize}
\item
$A_0\rta A_1\rta\cdots\rta A_n$ est une suite d'isog\'enies de sch\'emas
ab\'eliens de dimension $n$ sur $S=\Spec(R)$ tels que les 
${\rm Ker}(A_i\rta A_{i+1})$ sont des sch\'emas en groupes finis et
plats  de rang $p$ sur $S$
\item
$\lambda_0$ et $\lambda_n$ sont des polarisations principales de $A_0$ et
$A_n$ telles que le compos\'e des fl\`eches partant de $A_i$ et revenant
\`a $A_i$ apr\`es avoir fait le tour du diagramme 
$$\begin{array}{ccccccc}
A_0&{\hfl\alpha}&A_{1}&{\hfl\alpha}&\cdots&{\hfl\alpha}&A_n\\
\vflm{\lambda_0^\smv}{}&&&&&&\vfld{}{\lambda_n}\\ 
A^\smv_0&\hflg{}{\alpha^\smv}&A^\smv_{1}&\hflg{}{\alpha^\smv}
&\cdots&\hflg{}{\alpha^\smv}&A^\smv_{n}
\end{array}$$ est \'egal \`a $p.{\rm Id}_{A_i}$ pour tout
$i=0,\ldots,g$. Ici on a d\'egign\'e par $A^\smv_i$ le sch\'ema
ab\'elien dual \`a $A_i$.
\item $\iota_N$ est un isomorphisme symplectique $\iota_N: A_0[N]\rta
(\ZZ/N\ZZ)^{2n}$ pour une forme symplectique fix\'ee sur
$(\ZZ/N\ZZ)^{2n}$ 
\end{itemize}

Pour $N$ assez grand, ce foncteur est repr\'esentable par un
$T$-sch\'ema ${{\cal A}}$ ayant des singularit\'es dans la fibre
sp\'eciale. Rappelons la construction de de Jong  et de Rapoport-Zink du
mod\`ele local de ces singularit\'es.
 
Pour tout $i=0,\ldots,n$, les ingr\'edients principaux de la construction
du mod\`ele local sont $M_i= \rmR^1 {a_i}_*(\Omega^\bullet_{A_i/S})$ et
$\omega_i =a_{i,*}\Omega^1_{A_i/S}$ o\`u $a_i$ est le morphisme
structurel
$a_i:A_i\rta S$.  

Les $M_i$ sont des $\calO_S$-module localement libre de rang $2n$ qui
viennent avec :
\begin{itemize}
\item les homomorphismes de $\calO_S$-modules
$$M_n\hfl\alpha M_{n-1}\hfl\alpha\cdots\hfl\alpha M_0$$ induits des
isog\'enies $\alpha:A_i\rta A_{i+1}$.  

\item les formes symplectiques non d\'eg\'en\'er\'ees
$$\begin{array}{ccc} q_0:M_0\otimes_{\calO_S} M_0\rta \calO_S &\mbox{et}
&  q_n:M_n\otimes_{\calO_S}M_n\rta \calO_S
\end{array}$$ induites par les polarisations principales $\lambda_0$ et
$\lambda_n$ sur $A_0$ et $A_n$. 
\end{itemize} Notons $M(A)$ la donn\'ee
$$M(A)=(M_n\hflg\alpha{} M_{g-1}\hflg\alpha{}\cdots\hflg\alpha{}
M_0,q_0,q_n ).$$ Elle v\'erifie les conditions suivantes (voir loc. cit.
prop. 3.1)
\begin{itemize}
\item Les ${\rm Coker}(M_i\hfl\delta M_{i+1})$ sont des
$\calO_S/p\calO_S$-module localement libres de rang $1$.
\item Pour tout  $i=0,\ldots,n$, le compos\'e des fl\`eches partant de
$M_i$  et revenant \`a $M_i$ apr\`es avoir fait le tour du diagramme 
$$\begin{array}{ccccccc}
M_0&\hflg{\alpha}{}&M_{1}&\hflg\alpha{}&\cdots&\hflg\alpha{}&M_n\\
\vfld{q_0^\smv}{}&&&&&&\vflm{}{q_n}\\ 
M^\smv_0&\hfl{\alpha^\smv}&M^\smv_{1}&\hfl{\alpha^\smv}
&\cdots&\hfl{\alpha^\smv}&M^\smv_{n}
\end{array}$$ est \'egal \`a la multiplication par $p.{\rm Id}_{M_i}$.
Ici on a pos\'e
$M_i^\smv={\rm Hom}(M_i,\calO_S)$. 
\end{itemize}

On appelle la donn\'ee $M$ de $\calO_S$-modules localement libres de
rang $2n$   munis des homomorphismes $\calO_S$-lin\'eaires $\alpha$ 
$$(M_0\hflg\alpha{} M_{1}\hflg\alpha{}\cdots\hflg\alpha{} M_n,q_0,q_n)$$
et des formes symplectiques $q_0$ et $q_n$ v\'erifiant les deux
derni\`eres propri\'et\'es, {\em un syst\`eme {\rm Sp} de
$\calO_S$-modules de type
$\Gamma_0(p)$}. Rappelons le lemme suivant, d\^u \`a de Jong
(loc. cit.  lemma  3.6 ; voir aussi
\cite{RZ}, {\it appendix to chapter 3: normal forms of lattice chains}
pour un \'enonc\'e pus g\'en\'eral).

\begin{lemme}[de Jong] Deux syst\`emes {\rm Sp} de $\calO_S$-modules de
type $\Gamma_0(p)$ arbitraires $M$ et $M'$ sont localement isomorphes
pour la topologie de Zariski de $S$.  Si $I$ est un id\'eal de carr\'e
nul de $\calO_S$, si $S'$ est le sous-sch\'ema ferm\'e de $S$ d\'efini
par $I$, tout isomorphisme entre les restrictions de $M$ et $M'$ \`a $S'$
se rel\`eve en un isomorphisme entre $M$ et $M'$.
\end{lemme}

Il existe un syst\`eme Sp de $\calO_T$-modules de type $\Gamma_0(p)$
standard, not\'e $V$ dont on pr\'ecisera la construction dans la section
suivante. Une cons\'equence imm\'ediate du lemme de de Jong est que le
foncteur $S\mapsto \Aut(V\otimes_{\calO_T} \calO_S)$  est
repr\'esentable par un sch\'ema en groupes lisse de type fini sur $T$. 

On consid\`ere le foncteur $\calW$ qui associe \`a tout $T$-sch\'ema $S$
l'ensemble  des couples $(A,\iota)$ o\`u $A\in {{\cal A}}(S)$ et o\`u
$\iota$ est un isomorphisme
$$\iota: M(A)\rta V\otimes_{\calO_T}\calO_S$$ de syst\`eme Sp de
$\calO_S$-modules de type $\Gamma_0(p)$.  Le foncteur d'oubli induit un
morphisme $\pi:\calW\rta {{\cal A}}$. L'\'enonc\'e suivant r\'esulte
alors \'egalement du lemme de de Jong rappel\'e ci-dessus.

\begin{proposition} Le morphisme $\pi:\calW\rta {{\cal A}}$ est un
torseur sous le sch\'ema en groupes
$\Aut(V)$. En particulier, c'est un morphisme repr\'esentable,
s\'epar\'e lisse et  surjectif. 
\end{proposition}

Maintenant pour chaque $i$, $M_i=\rmR^1 {a_i}_* \Omega^\bullet _{A_i/S}$
contient un 
$\calO_S$-sous-module 
$$\omega_i= {a_i}_* \Omega_{A_i/S}^1\subset M_i$$ qui est localement un
facteur direct de rang $n$. Ces sous-modules $\omega_i$  v\'erifient les
conditions suivantes
\begin{itemize}
\item $\alpha(\omega_{i+1})\subset \omega_i$,
\item $\omega_0$ et $\omega_n$ sont totalement isotropes par rapport \`a 
$q_0$ et $q_n$ respectivement.
\end{itemize} Consid\'erons le foncteur $\calM$ qui associe \`a chaque
$T$-sch\'ema $S$ l'ensemble  des donn\'ees $L$ d'un
$\calO_S$-sous-modules $L_i\subset V_i$ localement facteur direct de
rang $n$ pour $i=0,\ldots,n$ tels que
\begin{itemize}
\item $\alpha(L_{i+1})\subset L_i$
\item $L_0$ et $L_n$ sont totalement isotropes par rapport \`a $q_0$ et
$q_n$ respectivement.
\end{itemize}

La donn\'ee d'un $S$-point $A\in {{\cal A}}(S)$ et d'un isomorphisme
$\iota:M(A)\rta V\otimes \calO_S$  de syst\`emes Sp de $\calO_S$-modules
de type $\Gamma_0(p)$, d\'efinit un point $L\in\calM(S)$  par la r\`egle
$L_i=\iota(\omega_i)$. On obtient ainsi un morphisme $f:\calW\rta \calM$.

En utilisant la  th\'eorie de Grothendieck-Messing des d\'eformations de
sch\'emas  ab\'eliens, de Jong \cite{deJong} et Rapoport et Zink
\cite{RZ} ont d\'emontr\'e le th\'eor\`eme suivant.

\begin{theoreme}[de Jong, Rapoport-Zink] Le morphisme $f:\calW\rta
\calM$ est un morphisme  lisse.
\end{theoreme}

D'apr\`es  \cite{Gen}, on sait de plus que le morphisme
$f:\calW\rta\calM$ est  {\em surjectif}.

\section{Syst\`eme standard} Soit $\bbA^1_T=\Spec(\ZZ_p[t])$ la droite
affine au-dessus de $T=\Spec(\ZZ_p)$. Soient 
$\calV_0,\ldots,\calV_{2n}$ des $\calO_{\bbA^1_T}$-modules libres de
rang $2n$. Soient 
$q_0:\calV_0\otimes\calV_0\rta\calO_{\bbA^1_T}$ et
$q_n:\calV_n\otimes\calV_n\rta
\calO_{\bbA^1_T}$ les formes symplectiques non-d\'eg\'en\'er\'ees
associ\'ees \`a la matrice  altern\'ee 
$$J=\left(\begin{array}{cc} 0 & K_n \\ -K_n & 0 \end{array}\right)$$
o\`u $K_n$ est la matrice $n\times n$ avec que des $1$ sur
l'anti-diagonale et que  des $0$ ailleurs. Soit $\alpha:\calV_{i+1}\rta
\calV_i$ l'homomorphisme de $\calO_{\bbA^1_T}$-modules d\'efini par 
$$\alpha=\left(\begin{array}{cc} 0 & \Id_{2n-1} \\ t-p & 0
\end{array}\right).$$

Notons $x_0:T\rta \bbA^1_T$ et $x_p:T\rta \bbA^1_T$ les sections
d\'efinies par $t=0$ et 
$t=p$ respectivement. Puisque $\alpha^{2n}=(t-p)\Id_{2n}$, les conoyaux
de $\alpha:
\calV_{i+1}\rta \calV_i$ sont support\'ees par la section $x_p$.

Posons maintenant $V_i= x_0^* \calV_i$ et d\'esignons aussi par
$\alpha,q_0,q_n$ les restrictions  de $\alpha,q_0,q_n$ \`a $x_0$. On
v\'erifie sans peine l'assertion suivante.

\begin{proposition} La donn\'ee 
$V=(V_0\hflg\alpha{} V_{1}\hflg\alpha{}\cdots\hflg\alpha{} V_n,q_0,q_n)$
forme un syst\`eme {\rm Sp} de $\calO_T$-modules de type $\Gamma_0(p)$.
\end{proposition}

Il sera commode d'identifier $\calV_i$ avec son image par
l'homomorphisme injectif
$\alpha^i:\calV_i\rta \calV_0$ de sorte que la suite des morphismes 
$$\calV_{0} \hflg{\alpha}{} \calV_{1} \hflg{\alpha}{} \cdots
\hflg{\alpha}{} \calV_{2n}$$ puisse s'identifier \`a une suite
d'inclusions de $\calO_{\bbA^1_T}$-modules
$$\calV_{0}\supset \calV_{1} \supset \cdots \supset \calV_{2n-1}\supset
\calV_{2n}=\calV_0[-x_p] .$$ Apr\`es cette identification, la forme
$(t-p)^{-1}K$ d\'efinit des accouplements parfaits 
$\calV_i\otimes_{\calO_{\bbA^1_T}} \calV_{2n-i}\rta \calO_{\bbA^1_T}$
qui induisent des accouplements  parfaits $V_i\otimes_{\calO_T} V_{2n-i}
\rta \calO_T$. De plus, l'inclusion $\calV_i\subset \calV_{i-1}$ est
duale \`a $\calV_{2n-i+1}
\subset \calV_{2n-i}$. On est maintenant en mesure de  r\'e\'ecrire le
probl\`eme de modules de $\calM$  en termes de cha\^\i nes de r\'eseaux. 

Rappelons que le foncteur mod\`ele local $\calM$ associe \`a chaque
$T$-sch\'ema $S$ l'ensemble  des donn\'ees $L$ d'un
$\calO_S$-sous-modules $L_i\subset V_i\otimes_{\calO_T}\calO_S$
localement facteur direct de rang $n$ pour $i=0,\ldots,n$ tels que
\begin{itemize}
\item $\alpha(L_{i+1})\subset L_i$,
\item $L_0$ et $L_n$ sont totalement isotropes par rapport \`a $q_0$ et
$q_n$ respectivement.
\end{itemize} En utilisant la dualit\'e entre $V_i$ et $V_{2n-i}$, on
obtient pour tout $i=n+1,\ldots,2n$ un sous-$\calO_S$-module
$L_{i}\subset V_{i}$ qui est dual \`a $L_{2n-i}\subset V_{2n-i}$.

Notons que la donn\'ee d'un sous-$\calO_S$-module localement facteur
direct $L_i\subset  V_i\otimes_{\calO_T}\calO_S$ est \'equivalente \`a
la donn\'ee d'un $\calO_{\bbA^1_S}$-module localement libre de rang $2n$
munie des modifications
$$\calV_i\otimes_{\calO_T}\calO_S [-x_0] \subset \calL_i \subset
\calV_i\otimes_{\calO_T}\calO_S$$ si bien que la donn\'ee d'un $S$-point
de $\calM$ est maintenant \'equivalente \`a celle d'un  diagramme
d'inclusions de $\calO_{\bbA^1_S}$-modules 
$$\begin{array}{rcccccl}
\calV_{0,S}& \supset& \calV_{1,S}& \supset & \cdots& \supset
&\calV_{2n,S}=\calV_{0,S}[-x_p] \\
\cup & & \cup &&&& \cup\\
\calL_{0}& \supset& \calL_{1}& \supset & \cdots& \supset
&\calL_{2n}=\calL_0[-x_p]\\
\cup & & \cup &&&& \cup \\
\calV_{0,S}[-x_0]& \supset& \calV_{1,S}[-x_0]& \supset & \cdots& \supset
&\calV_{2n,S}[-x_0] \\
\end{array}$$ telle que pour tout $i$, les deux modules $\calL_i$ et
$\calL_{2n-i}$ sont en dualit\'e par  rapport \`a la forme $(t-p)^{-1}K$
et que $\calL_i/\calV_{i.S}[-x_0]$ est un sous-$\calO_S$-module 
localement facteur direct de rang $n$ de
$V_{i,S}=\calV_{i,S}/\calV_{i,S}[-x_0]$. 

Cette nouvelle pr\'esentation du probl\`eme de modules $\calM$ semble
nettement plus compliqu\'ee que la pr\'ec\'edente. En contre partie, elle
met en 
lumi\`ere le fait important suivant (pour lequel nous renvoyons \`a
l'article de Kottwitz et Rapoport
\cite{KR}, \`a celui de Goertz
\cite{Goertz} et \`a celui de Haines et du premier auteur 
\cite{HN} ; voir aussi dans \cite{Gen} la preuve de la surjectivit\'e du
morphisme
$f:\calW \rta \calM$). 

\begin{proposition} La fibre sp\'eciale $\calM_s$ est un sous-sch\'ema
ferm\'e  de l'ind-sch\'ema  de drapeaux affine de
${\rm GSp}(2n)$ (au sens de \cite{BL}).
\end{proposition}

Pour tout couple d'entiers $r_-<r_+$, consid\'erons le ${{\mathbb
F}}_p$-sch\'ema $X_{r_{\pm}}$ dont  l'ensemble des $S$-points,  pour
toutes ${{\mathbb F}}_p$-alg\`ebre $S$, est l'ensemble des diagrammes
d'inclusions de $\calO_{\bbA^1_S}$-modules 
$$\begin{array}{rcccccl}
\calV_{0,S}[r_+x]& \supset& \calV_{1,S}[r_+x]& \supset & \cdots& \supset 
&\calV_{2n,S}[r_+x] \\
\cup & & \cup &&&& \cup\\
\calL_{0}& \supset& \calL_{1}& \supset & \cdots& \supset
&\calL_{2n}=\calL_0[-x]\\
\cup & & \cup &&&& \cup \\
\calV_{0,S}[r_-x]& \supset& \calV_{1,S}[r_-x]& \supset & \cdots& \supset
&\calV_{2n,S}[r_-x] \\
\end{array}$$  o\`u $x=x_0=x_p$ en caract\'eristique $p$, et tels que
$\calL_i$ et
$\calL_{2n-i}$ sont duaux par rapport \`a la forme $t^{r_++r_-}K$ et que
$\calL_i/\calV_{i,S}[r_-x]$ est un sous-$\calO_S$-module localement
facteur direct de rang $n(r_+-r_-)$ de
$\calV_{i,S}[r_+x]/ \calV_{i,S}[r_-x]$. Visiblement $\calM_s=X_{0,-1}$.
Les 
$X_{r_\pm}$ sont tous des ${{\mathbb F}}_p$-sch\'emas projectifs et
forment un syst\`eme inductif avec $r_+\rta \infty$ et $r_-\rta
-\infty$.  

Si on localise en dehors de la section $x$, les inclusions dans le
diagramme deviennent toutes des \'egalit\'es. Notons $\calO_x$ l'anneau
local compl\'et\'e de ${{\mathbb F}}_p[t]$ en $t=0$. On peut identifier
les compl\'et\'es en $x$ de $\calV_{i,{{\mathbb F}}_p}$ avec  
$$\calV_{i,{{\mathbb F}}_p,x}=\bigoplus_{j=1}^i\calO_x e_j\oplus
\bigoplus_{j=i+1}^{2n} t\calO_x e_j$$  o\`u $e_1,\ldots,e_{2n}$ est la
base standard. Par le th\'eor\`eme de recollement formel de Beauville et
Laszlo \cite{BL}, la donn\'ee d'un
${{\mathbb F}}_p$-point de $\calM_{r_\pm}$ est \'equivalente \`a la
donn\'ee d'une cha\^\i ne p\'eriodique de r\'eseaux 
$$\calL_{0,x} \supset \calL_{1,x} \supset  \cdots \supset \calL_{2n,x}=t
\calL_{0,x}$$  v\'erifiant la condition d'inclusion 
$$t^{-r_+}\calV_{i,{{\mathbb F}}_p,x}\supset \calL_{i,x} \supset
t^{-r_-}\calV_{i,{{\mathbb F}}_p,x}$$  pour tout $i$ et la condition de
dualit\'e \'evidente. En prenant la limite quand $r_\pm\rta \pm\infty$
on obtient donc  l'ind-sch\'ema $X$ des drapeaux affines de ${\rm
GSp}(2n)$  
$$X({{\mathbb F}}_p)={\rm GSp}(2n,{{\mathbb F}}_p(\!(t)\!))/I({{\mathbb
F}}_p)$$ o\`u $I$ est un groupe alg\'ebrique sur ${{\mathbb F}}_p$ dont
les ${{\mathbb F}}_p$-points forment le sous-groupe d'Iwahori (oppos\'e)
standard de ${\rm GSp}(2n,{{\mathbb F}}_p(\!(t)\!))$ --la mention
``oppos\'e" r\'esulte de notre choix de travailler avec les cha\^\i nes
descendantes. 

En particuler, le ind-groupe ${\rm GSp}(2n,{{\mathbb F}}_p(\!(t)\!))$
agit sur
$X$. L'action restreinte \`a $I$ laisse stable les $X_{r_\pm}$ et en
particulier, $I$ agit sur $X_{0,-1}=\calM_s$. Il existe en fait un
homomorphisme canonique de groupes alg\'ebriques 
$$I\rta {\rm Aut}(V_s)$$
\`a travers lequel se factorise l'action de $I$ sur $\calM_s$. Pour la
v\'erification facile, mais fastidieuse, des assertions pr\'ec\'edentes,
nous renvoyons \`a \cite{HN}. 

\section{L'ensemble de Kottwitz-Rapoport}

D'apr\`es Iwahori-Marsumoto \cite{IM}, on a la d\'ecomposition
$${\rm GSp}(2n,{{\mathbb F}}_p(\!(t)\!))=\bigsqcup_{w\in \widetilde W_a}
I({{\mathbb F}}_p)w I({{\mathbb F}}_p)$$  o\`u $\widetilde W_a$ est le
groupe de Weyl affine \'etendu du groupe ${\rm GSp}(2n)$. Cette
d\'ecomposition se traduit g\'eom\'etriquement en une d\'ecomposition
cellulaire $X=\bigsqcup_{w\in \widetilde W_a} X_w$ avec
$X_w=\bbA^{\ell(w)}$ o\`u $\ell(w)$ est la fonction longueur dont on
rappellera la d\'efinition juste apr\`es.  Du simple fait que $\calM_s$
est un sous-sch\'ema projectif et $I$-\'equivariant de $X$, on a 
$$\calM_s=\bigsqcup_{w\in\KR(\mu)} X_w$$  o\`u $\KR(\mu)$ est un certain
sous-ensemble fini de $\tilde W_a$ que Kottwitz et Rapoport appellent
l'ensemble des \'el\'ements
$\mu$-permissibles. Pour rappeler la d\'efinition combinatoire, due \`a
Kottwitz et Rapoport \cite{KR} de cet ensemble fini, il nous faut faire
quelques rappels sur le groupe de Weyl affine. 

Soit $G$ un groupe semi-simple d\'eploy\'e d\'efini sur $\ZZ_p$, $T$ un
tore maximal d\'eploy\'e de $T$, $W={\rm Nor}(T)/T$ le groupe de Weyl
fini associ\'e. Notons $Q$ le $\ZZ$-module engendr\'e par les racines de
$G$ dans l'espace vectoriel r\'eel $X_\RR={\rm Hom}(T,\GG_m)\otimes
\RR$.  Notons $Q^\smv$ le $\ZZ$-module engendr\'e par les coracines dans
l'espace vectoriel r\'eel dual $X^\smv_\RR={\rm Hom}(\GG_m,T)\otimes
\RR$ et
$P^\smv$ le r\'eseau dual \`a $Q$. 

Soient $W_a=W\ltimes Q^\smv$ le groupe de Weyl affine et $\widetilde
W_a=W\ltimes P^\smv$ le groupe de Weyl affine \'etendu. Rappelons que le
groupe $W_a$ est un groupe de Coxeter de sorte qu'il muni d'un ordre
partiel --celui de Bruhat-- et d'une fonction longueur $\ell:W_a \rta
\NN$. Une racine $\alpha$ de $G$ et un entier $k\in\ZZ$ d\'efinissent un
mur affine 
$$H_{\alpha,k}=\{x\in X^\smv_\RR\mid \langle\alpha,x\rangle =k\}.$$  Le
compl\'ementaire dans $X^\smv_\RR$ de la r\'eunion de tous les murs
affines se d\'ecompose en composantes connexes qui sont appel\'ees les
alcoves.  

Soient $(\alpha_i)_{i\in I}$ les racines simples associ\'ees au choix
d'un sous-groupe de Borel standard, $\gamma$ la plus grande racine. Il
existe une alcove, dite de base,  qui est d\'elimit\'ee par les murs
$H_{\alpha_i,0}$ pour $i\in I$ et $H_{\gamma,1}$.  Pour chaque mur de
l'alcove, il existe un unique sommet, dit oppos\'e, qui ne lui
appartient pas. Notons pour tout $i\in I$, $a_i$ le sommet oppos\'e
\`a $H_{\alpha_i,0}$ et notons $a_0$ le sommet oppos\'e \`a
$H_{\gamma,1}$.

Le groupe $\widetilde W_a$ agit sur $X^\smv_\RR$   en envoyant
un mur affine sur un autre, de sorte qu'il agit sur l'ensemble des
alcoves. Rappelons que l'action de $W_a$ sur l'ensemble des alcoves est
simple et transitive, voir \cite{Hum}. Il s'ensuit que $\widetilde
W_a=W_a
\Omega$ o\`u $\Omega$ est le fixateur de l' alcove de base qui est un
sous-groupe fini. L'ordre de Bruhat et la fonction longueur s'\'etendent
trivialement de $W_a$ \`a $\widetilde W_a$. Pour tous $w,w'\in W_a$,
$\tau,\tau'\in \Omega$ on dit $w\tau\leq w'\tau'$ si et seulement  si
$w\leq w'$ et $\tau=\tau'$. On pose $\ell(w\tau)=\ell(w)$.

Soit $\mu\in Q^\smv$ un copoids minuscule de $G$. D'apr\`es Kottwitz et
Rapoport, un \'el\'ement $w\in \widetilde W_a$ est dit $\mu$-{\em
permissible} si pour tout sommet $a$ de l'alcove de base, la diff\'erence
$w(a)-a$ est conjugu\'ee \`a $\mu$ par rapport \`a l'action de groupe de
Weyl fini $W$. 

Pour comprendre g\'eom\'etriquement cette condition, consid\'erons le cas
${\rm GL}(d)$. Le groupe de Weyl affine \'etendu $\widetilde W_a({\rm
GL}(d))$ agit sur l'ensemble des alcoves de ${\rm PGL}(d)$ \`a travers
l'homomorphisme canonique surjectif $\widetilde W_a({\rm GL}(d)) 
\rta \widetilde W_a({\rm PGL}(d))$. Soit 
$$\mu=(\underbrace{1,\ldots,1}_r,0,\ldots,0)$$ un copoids minuscule de
${\rm GL}(d)$ et notons aussi par $\mu$ le copoids minuscule induit de
${\rm PGL}(d)$. On choisit pour chaque $w\in\widetilde W({\rm GL}(d))$
le repr\'esentant, not\'e aussi $w\in{\rm GL}(d,{{\mathbb
F}}_q(\!(t)\!)))$ la matrice $d\times d$ avec exactement $d$ entr\'ees
non nulles lesquelles sont de la forme $t^s$. Ce choix est
caract\'eris\'e par la propri\'et\'e que $w$ stabilise l'ensemble de
vecteurs $t^m e_j$ o\`u $\{e_j\}$ est la base standard.  

Rappelons qu'un ${{\mathbb F}}_p$-point $\calL$ de l'ind-sch\'ema des
drapeaux affines de ${\rm GL}(d)$ est un drapeau p\'eriodique de
r\'eseaux
$$\calL_{0,x}\supset \calL_{1,x}\supset\cdots\supset
\calL_{d,x}=t\calL_{0,x}$$  dans ${{\mathbb F}}_p(\!(t)\!)^{d}$.
Rappelons aussi qu'on a d\'efini un drapeau standard $\calV_{i,{{\mathbb
F}}_p,x}$ dans la section 2.

\begin{proposition} Supposons que $\calL\in X_w(k)$ pour un certain
$w\in \widetilde W_a({\rm GL}(d))$ avec ${\rm val}(\det(w))=r$. Alors,
$w$ est $\mu$-permissible si et seulement si pour tout $i=0,\ldots,d-1$,
on a $\calV_{i,{{\mathbb F}}_p,x}\supset
\calL_i\supset t\calV_{i,{{\mathbb F}}_p,x}$. 
\end{proposition}

\dem Par d\'efinition, $I$ stabilise le drapeau de base
$\calV_{i,{{\mathbb F}}_p,x}$, si bien qu'il suffit de tester le cas
$\calL_{i,x}=w \calV_{i,k,x}$. On doit  d\'emontrer que pour tout
$i=0,\ldots,d-1$, les inclusions
$$\calV_{i,{{\mathbb F}}_p,x}\supset w\calV_{i,x} \supset
t\calV_{i,{{\mathbb F}}_p,x}$$ ont lieu si et seulement si
$w(a_i)-a_i\in W\mu$, toujours sous l'hypoth\`ese
${\rm val}(\det(w))=r$. 

Du fait qu'on consid\`ere le cas lin\'eaire, les sommets
$a_1,\ldots,a_{d-1}$ de l'alcove de base sont les copoids fondamentaux.
Notons ${\rm t}_{a_i}$ l'\'el\'ement $a_i$ vu comment un
\'el\'ement de $\widetilde W_a$ et comme \'el\'ement de ${\rm
GL}(d,k(\!(t)\!))$. On a $\calV_{i,{{\mathbb F}}_p,x}={\rm
t}_{a_i}\calV_{0,{{\mathbb F}}_p,x}$. Les inclusions  
$$\calV_{i,{{\mathbb F}}_p,x}\supset w\calV_{i,x} \supset
t\calV_{i,{{\mathbb F}}_p,x}$$ ont lieu si et seulement si 
$${\rm t}_{-a_i}w{\rm t}_{a_i}\in K\mu K$$ o\`u $K={\rm GL}(n,k[[t]])$
est le stabilisateur de $\calV_{0,k}$. Il s'ensuit que ${\rm
t}_{-a_i}w{\rm t}_{a_i}\in y{\rm t}_\mu W$ pour un certain $\mu\in W$
d'o\`u $w(a_i)-a_i=y\mu$ en appliquant ${\rm t}_{-a_i}w{\rm t}_{a_i}$
\`a $0$. \findem  

Notons tout de suite que cette d\'emonstration n'est valable que dans le
cas lin\'eaire. Dans le cas du groupe ${\rm GSp}(2n)$ avec le copoids
minuscule  
$$\mu=(\underbrace{1,\ldots,1}_n,\underbrace{0,\ldots,0}_n)$$ on
d\'eduit directement de la proposition pr\'ec\'edente et de la
description fonctorielle de $\calM_s$ l'assertion suivante. 

\begin{corollaire} On a une d\'ecomposition cellulaire
$$\calM_s=\bigsqcup_{w\in \KR(\mu)} X_w$$ o\`u $\KR(\mu)$ est l'ensemble
des \'el\'ements de $\widetilde W_a({\rm GSp}(2n))$ dont l'image dans
${\widetilde W}_a({\rm GL}(2n))$ est $\mu$-permissible. 
\end{corollaire}

Kottwitz et Rapoport ont d\'emontr\'e qu'en plus, $\KR(\mu)$ est
exactement l'ensemble des $\mu$-permissibles dans ${\rm GSp}(2n)$ car
les
\'el\'ements qui sont $\mu$-permissibles dans ${\rm GSp}(2n)$ sont aussi
ceux dont l'image est $\mu$-permissible dans ${\rm GL}(2n)$.  Nous
pr\'ef\'erons la description du corollaire qui est moins \'el\'egant,
mais qui colle automatiquement \`a la description fonctorielle de
$\calM_s$. Par ailleurs, comme cela a \'et\'e montr\'e dans \cite{HN2},
pour un $\mu$ non minuscule, il existe des \'el\'ements
$\mu$-permissibles dans ${\rm GSp}(2n)$ dont l'image n`est pas
$\mu$-permissible dans ${\rm GL}(2n)$. Mais revenons au cas minuscule
pour rappeller le r\'esultat important suivant.

\begin{theoreme}[Kottwitz-Rapoport] Si $w\in\KR(\mu)$ alors il existe
$\rm y\in W({\rm GSp}(2n))$ --le groupe de Weyl vectoriel de ${\rm
GSp}(2n)$-- tel que $w\leq {\rm t}_{\rm y\mu}$. 
\end{theoreme}

Les \'el\'ements translations ${\rm t}_{{\rm y}\mu}$ sont clairement
$\mu$-permissibles parce que ${\rm t}_{{\rm y}\mu}(a_i)-a_i={\rm y}\mu$.
Le th\'eor\`eme dit qu'ils  sont exactement les \'el\'ements maximaux de
$\KR(\mu)$. Cet \'enonc\'e a
\'et\'e g\'en\'eralis\'e par Haines et Ng\^o \`a tous les poids dominants
de ${\rm GSp}(2n)$ \cite{HN2}.

\section{Alcoves et $p$-rang}

Revenons aux notations de la section 1 o\`u on a construit un morphisme
lisse ${\rm Aut}(V)$-\'equivariant $f:\calW\rta \calM$ et aussi un
morphisme $\pi:\calW \rta {{\cal A}}$ qui fait de $\calW$ un ${\rm
Aut}(V)$-torseur au-dessus de ${{\cal A}}$. La stratification en ${\rm
Aut}(V)_s$-orbites de $\calM_s$ d\'efinie dans la section 3
$$\calM_s=\bigsqcup_{w\in\KR(\mu)} X_w$$ induit par image inverse une
stratification $\calW_s=\bigsqcup_{w\in
\KR(\mu)} \calW_w$ dont les strates $\calW_w=f^{-1}(X_w)$ sont ${\rm
Aut}(V)_s$-\'equivariantes. Du fait que $\calW$ est un ${\rm
Aut}(V)$-torseur au-dessus de ${{\cal A}}$, cette stratification se
descend en une stratification de ${{\cal A}}_s$ 
$${{\cal A}}_s=\bigsqcup_{w\in\KR(\mu)}{{\cal A}}_w$$ telle que
$\calW_w=\pi^{-1}({{\cal A}}_w)$ pour tout $w\in \KR(\mu)$. Cette
stratification jouit des propri\'et\'es tr\`es agr\'eables, d\'eduites de
celles connues sur le mod\`ele local, que nous \'enum\'erons ici pour des
r\'ef\'erences ult\'erieures :  
\begin{itemize}
\item Les strates ${{\cal A}}_w$ sont lisses de dimension $l(w)$
d'apr\`es le th\'eor\`eme de lissit\'e de de Jong et Rapoport-Zink.
Elles sont non vides d'apr\`es la surjectivit\'e de (\cite{Gen},
proposition 1.3.2).
\item La restriction du complexe de cycles proches $\rmR\Psi_{{\cal A}}
\QQ_\ell$ \`a chaque strate ${{\cal A}}_w$ est constante. En effet, du
c\^ot\'e du mod\`ele local $\calM$, le complexe des cycles proches est
constant sur les strates $\calM_w$ parce qu'il est \'equivariant par
rapport \`a l'action de ${\rm Aut}(V)$.  
\item La trace semi-simple de Frobenius agissant sur une fibre quelconque
de $\rmR\Psi_{{\cal A}} \QQ_\ell |_{{{\cal A}}_w}$ est calcul\'ee par la
fonction de Bernstein-Kottwitz, d'apr\`es Haines et Ng\^o \cite{HN}. On
a par ailleurs une formule tr\`es explicite de cette fonction due \`a
Haines \cite{Haines} : pour tout $A\in {{\cal A}}_w({{\mathbb F}}_q)$   
$${\rm Tr}^{ss}\left({\rm Frob}_{q^n},(\rmR\Psi_{{\cal A}}
\QQ_\ell)_A\right) =\epsilon_{t(w)}\epsilon_{w} R_{w,{\rm
t}_\lambda}(q^n)$$ o\`u $R_{w,{\rm t}(w)}$ est le polyn\^ome $R_{x,y}$
habituel de Kazhdan-Lusztig, o\`u $\lambda$ est l'unique \'el\'ement de
$P^\smv$ vu comme sous-groupe de $\widetilde W$ tel que $w=x{\rm
t}_\lambda$  avec $x\in W$ un \'el\'ement de groupe de Weyl vectoriel et
o\`u 
$\epsilon_w$  est le signe $(-1)^{\ell(w)}$. 
\end{itemize} Il est naturel d'essayer d'interpr\'eter directement cette
stratification en termes de vari\'et\'es ab\'eliennes. 

Soit $A$ une vari\'et\'e ab\'elienne de dimension $n$ d\'efinie sur un
corps s\'eparablement clos $k$ de caract\'eristique $p$. Par
d\'efinition,  le {\em $p$-rang} de $A$ est l'entier $r$ tel que
$\#A[p](k)=p^r$ o\`u
$A[p]$ est le sous-groupes des $p$-torsions de $A$. Le $p$-rang est un
entier compris entre $0$ et $n$. La vari\'et\'e ab\'elienne $A$  est
dite {\em ordinaire} si son $p$-rang vaut $n$. Le $p$-rang ne  d\'epend
de $A$ qu'\`a l'isog\'enie pr\`es si bien  qu'il d\'efinit une
fonction sur l'ensemble des points g\'eom\'etriques de ${{\cal A}}_s$.

Soit $w\in \KR(\mu)$ vu comme un \'el\'ement de groupe de Weyl affine
\'etendu. On peut le repr\'esenter comme une matrice $2n\times 2n$
normalisant le tore diagonal et ayant donc exactement $2n$ entr\'ees non
nulles,
$n$ d'entre elles  valant $1$ et les $n$ autres valant $t$. On
d\'efinit $r(w)$ comme le nombre de $t$ sur la diagonale.

Une autre fa\c con agr\'eable de d\'efinir $r(w)$ est la suivante.  Son
image $x$ dans le groupe de Weyl vectoriel de ${\rm GSp}(2n)$ peut 
\^etre consid\'er\'e comme une permutation de l'ensemble
$\{1,\ldots,2n\}$ commutant avec l'involution $(2n,2n-1,\cdots,2,1)$.
L'ensemble des points fixes de $x$ est en particulier stable par cette
involution. Son cardinal est donc un nombre pair (puisque l'involution
est sans  point fixe)  et vaut en fait
$2r(w)$.

\begin{theoreme} La fonction $p$-rang est constante sur chaque strate
${{\cal A}}_w$ pour tout
$w\in \KR(\mu)$. Plus pr\'ecis\'ement, le $p$-rang d'un point
$A\in{{\cal A}}_w(k)$ est \'egal \`a l'entier $r(w)$.
\end{theoreme}

\dem Soit $A\in {{\cal A}}(k)$ o\`u $k$ est le corps s\'eparablement
clos de caract\'eristique $p$ figurant dans l'\'enonc\'e
$$A=\left(A_0\hfl\alpha A_1\hfl\alpha\cdots\hfl\alpha A_n
,\lambda_0,\lambda_n,\iota_N\right).$$ L'isog\'enie $A_{i-1}\hfl\alpha
A_i$ induit un morphisme lin\'eaire de rang
$2n-1$ entre les $k$-espaces $\alpha: M_i\rta M_{i-1}$ qui sont les
$\rmH^1_{\rm DR}$ de $A_i$ et $A_{i-1}$ respectivement. Ce morphisme
lin\'eaire induit par restriction un morphisme entre les sous-espaces
vectoriels $\omega_i\rta \omega_{i-1}$ de $M_i$ et $M_{i-1}$
respectivement. Soit $H_i$ le noyau de l'homormorphisme
$\alpha:A_{i-1}\rta A_i$. Ce $H_i$ \'etant un groupe fini et plat de
rang $p$, on peut se demander s'il est \'etale ou de type multiplicatif
ou biconnexe. D'apr\`es de Jong
\cite{deJong}, on sait que
\begin{enumerate}
\item $H_i$  est \'etale si et seulement si $\omega_i \rta \omega_{i-1}$
est un isomorphisme, ou ce qui est \'equivalent $\omega_i\cap{\rm ker} 
(M_i\rta M_{i-1})=0$,
\item $H_i$ est de type multiplicatif si et seulement si $\omega_{i-1}
\not\subset\alpha(M_i)$.
\end{enumerate} Ces deux conditions  s'\'echangent clairement par
dualit\'e,
$\omega_{i-1}^\smv$ a une intersection nulle avec ${\rm
ker}(\alpha^\smv:M^\smv_{i-1}\rta M^\smv_i)$ si et seulement si
$\omega_{i-1}\not\subset\alpha(M_i)$. La premi\`ere assertion vient du
fait que l'espace des diff\'erentielles sur $H_i$ est
$\omega_{i-1}/\alpha(\omega_i)$ et qu'il est nul si et seulement si $H_i$
est  \'etale.

Soit $(A,\iota)\in \calW(k)$ un rel\`evement de $A$ \`a $\calW$. Soit
$\calL=f({{\cal A}},\iota)\in\calM(k)$ dont la donn\'ee est
\'equivalente \`a un drapeau p\'eriodique de r\'eseaux $\calL_{i,x}$ tel
que
$$\begin{array}{rcccccl}
\calV_{0,k,x}& \supset& \calV_{1,k,x}& \supset & \cdots& \supset
&\calV_{2n,k,x} \\ \cup & & \cup &&&& \cup\\ \calL_{0,x}& \supset&
\calL_{1,x}& \supset & \cdots& \supset &\calL_{2n,x}=t\calL_{0,x}\\ \cup
& & \cup &&&& \cup \\ t\calV_{0,k,x}& \supset& t\calV_{1,k,x}& \supset &
\cdots&  \supset &t\calV_{2n,k,x} \\
\end{array}$$ avec une condition de dualit\'e entre $\calL_{i,x}$ et
$\calL_{2n-i,x}$. Les conditions 1 et 2  pr\'ec\'edents sont alors
\'equivalentes respectivement aux conditions suivantes
\begin{enumerate}
\item $\calL_{i,x}/t\calV_{i,k,x}\rta \calV_{i-1,k,x}/t\calV_{i-1,k,x}$
est injectif,
\item $\calL_{i-1,x}\not\subset \calV_{i,k,x}$.
\end{enumerate} Ces conditions \'etant invariantes par rapport \`a
l'action de $I$,  le type (\'etale, multiplicatif ou biconnexe) de $H_i$
ne d\'epend donc que de l'orbite  de $I$. Il s'ensuit que le $p$-rang de
$A$ est constant sur chaque strate ${{\cal A}}_w$. 

Soient maintenant $\calL_{i,x}=w\calV_{i,k,v}$ pour tout $i=0,\ldots,2n$
avec $w\in\KR(\mu)$. Les rel\`evements de $w$ dans ${\rm GL}(2n,{{\mathbb
F}}_p(\!(t)\!))$ ont \'et\'e choisis de sorte que $w$ stabilise
l'ensemble des vecteurs de la forme $\{t^me_j\}$. Les cha\^\i nes de
r\'eseaux $(\calL_{i,x}=w\calV_{i,k,v})_i$ jouissent donc de la
propri\'et\'e suivante :   pour tout $i=1,\ldots,n$, il existe un unique
vecteur parmi les ceux de la forme $t^m e_j$, qui appartient \`a
$\calL_{i-1,x}$, mais qui n'appartient pas \`a $\calL_{i,x}$ ; de plus ce
vecteur est $w(e_i)$. Cette observation permet de d\'emontrer que {\em
la condition $1$ est \'equivalente \`a $w(e_i)=te_i$ ; la condition $2$
est \'equivalente \`a $w(e_i)=e_i$}. En effet :

\begin{enumerate}
\item Le noyau du morphisme $\calV_{i,k,x}/t\calV_{i,k,x}\rta
\calV_{i-1,k,x}/t\calV_{i-1,k,x}$ est un $k$-espace vectoriel de
dimension
$1$ engendr\'e par l'image du vecteur $te_i$. L'application
$\calL_{i,x}/t\calV_{i,k,x}\rta \calV_{i-1,k,x}/t\calV_{i-1,k,x}$ est
injective si et seulement si $\calL_{i,x}$ ne contient pas ce vecteur.
Mais $te_i$ appartient \`a $t\calV_{i-1,k,x}\subset
\calL_{i-1,x}$ ;  si de plus il n'appartient pas \`a $\calL_{i,x}$, il
doit
\^etre \'egal \`a $w(e_i)$. Donc, l'injectivit\'e de
$\calL_{i,x}/t\calV_{i,k,x}\rta \calV_{i-1,k,x}/t\calV_{i-1,k,x}$  est
\'equivalente \`a $w(e_i)=te_i$. 

\item Puisque $\calL_{i-1,x}\subset \calV_{i-1,k,x}$, la condition
$\calL_{i-1,x}\not\subset \calV_{i,k,x}$ est satisfaite si et seulement
si 
$e_i\in \calL_{i-1,x}$. Mais $e_i\notin \calV_{i,k,x}$ d'o\`u $e_i\notin 
\calL_{i,x}$. Donc $e_i\in \calL_{i-1,x}$ si et seulement si
$w(e_i)=e_i$. 
\end{enumerate}

Le $p$-rang \'etant au nombre des indices $i=1,\ldots,2n$ pour
lesquelles la condition 1 est satisfaite, il vaut $r(w)$ le nombre de
$t$ sur la diagonale dans l'expression matricielle de $w$. \findem

\begin{corollaire} Soit $A\in {{\cal A}}_{w}(k)$.  Alors $A$ est
ordinaire si et seulement si $w=\rm t_{y\mu}$ pour un certain $y\in
W({\rm GSp}(2n))$.
\end{corollaire}

\dem Soit $w={\rm t}_{\lambda}x$ avec $x\in W({\rm GSp}(2n))$. Suivant
la seconde description du nombre $r(w)$, le $p$-rang de $A$ vaut
$n$ si et seulement si $x=1$, donc $r=t_\lambda$. La condition de
$\mu$-permissibilit\'e entra\^ine alors que $\lambda=y\mu$ pour un
certain
$y\in W({\rm GSp}(2n))$. 
\findem

\bigskip Notons ${{\cal A}}^{\rm ord}$ le lieu ordinaire de ${{\cal
A}}_s$ ; on a donc montr\'e que ${{\cal A}}^{\rm ord}$ est la r\'eunion
des strates correspondant aux translations ${{\cal A}}_{{\rm
t}_{y\mu}}$.   En conjonction avec le th\'eor\`eme de Kottwitz-Rapoport
cit\'e dans le paragraphe  pr\'ec\'edent, nous en tirons le corollaire
suivant.

\begin{corollaire} La partie ordinaire ${{\cal A}}^{\rm ord}$ est dense
dans ${{\cal A}}_s$.
\end{corollaire}

\begin{small}

\end{small}
\hbox to \hsize{
\vbox to 5cm{
\hbox {B. C. Ng\^o :}
\hbox{\tt ngo@math.univ-paris13.fr}
\hbox{CNRS, UMR 7539}
\hbox{Universit\'e Paris-Nord} 
\hbox{D\'epartement de math\'ematiques}
\hbox{av. J.-B. Cl\'ement} 
\hbox{93430 Villetaneuse}
\hbox{FRANCE }
\hbox{}
\hbox{et}
\hbox{}
\hbox{IHES}
\hbox{35 route de Chartres}
\hbox{91440 Bures sur Yvettes}
\hbox{FRANCE}
\vfill}
\hfill
\vbox to 5cm{
\hbox{A. Genestier : }
\hbox{\tt Alain.Genestier@math.u-psud.fr}
\hbox{CNRS, UMR 8628}
\hbox{Universit\'e Paris-Sud}
\hbox{D\'epartement de math\'ematiques}
\hbox{B\^atiment 425}
\hbox{91405 Orsay}
\hbox{FRANCE }
\vfill}
}

\end{document}